\documentclass{article}

\usepackage{amssymb, amsmath, amscd, amsthm, amsfonts}

\newcommand\Fp{\mathbb{F}_p}

\newcommand\Z{\mathbb{Z}}
\newcommand\Q{\mathbb{Q}}

\newcommand\F{{\mathcal{F}}}
\newcommand\A{{\mathcal{A}}}

\newcommand\R{{\mathbb{R}}}
\newcommand\Pp{{\mathcal{P}}}

\newcommand{\C}{{\mathbb{C}}}

\newcommand{\Qq}{\mathcal{Q}}

\newcommand{\Aa}{\alpha}

\newcommand{\Gg}{\Gamma}

\newcommand{\Oo}{\mathcal{O}}

\newcommand{\Wrr}{\Z/k\Z \text{ wr } S_n}

\newtheorem{theorem}{Theorem}[section]
\newtheorem{theorem'}{Theorem}
\newtheorem{lemma}[theorem]{Lemma}

\newtheorem{corollary}[theorem]{Corollary}

\newenvironment{example}{\refstepcounter{theorem}\smallskip
\noindent {\textbf{Example~\thetheorem}\ }}{\hskip\hsize plus0pt minus\hsize
\hbox{$\Box$}\bigskip}

\newenvironment{example'}{\smallskip
\noindent {\textbf{Example \ref{Galois} Continued~}\ }}{\hskip\hsize plus0pt minus\hsize
\hbox{$\Box$}\bigskip}

  \newenvironment{Proof}{\noindent{\bf Proof} \ }{\QED}\smallskip
\newcommand\QED{\newline \rightline{$\blacksquare$} \bigskip}

  \smallskip

\title{The Proportion of $k$-cycles for Polynomials Modulo Primes}

\author{Jonathan Root}
\date{}

\begin{document}

\maketitle






\begin{abstract}
Let $f(x) \in \Fp[x]$, and define the {\it orbit} of $x\in \Fp$ under the iteration
of $f$ to be the set
\[
\Oo(x):=\{x,f(x),(f\circ f)(x),(f\circ f\circ f)(x),\dots\}.
\]
An orbit is a {\it$k$-cycle} if it is periodic of length $k$.
In this paper we fix a polynomial $f(x)$ with integer coefficients
and for each prime $p$ we consider $f(x) \pmod p$ obtained by 
reducing the coefficients of  $f(x)$ modulo $p$.  We ask for the
density of primes $p$ such that $f(x)\pmod p$ has a $k$-cycle
in $\Fp$.  We prove that in many cases the density is at most
$1/k$.  We also give an infinite family of polynomials in each degree
with this property.

\end{abstract}







\section{Introduction}

If $K$ is a field and $f(x)\in K[x]$, 
then $\Aa \in K$ is a {\it periodic} point of $f$ if there exists an integer $k\geq 1$
such that $\Aa$ is a fixed point under the $k$-fold composition of $f$, which we denote by $f^k$.
A periodic point $\Aa$ has {\it minimal period} $k$ if $f^k(\Aa)=\Aa$,
but $f^j(\Aa) \neq \Aa$ for any $1\leq j<k$. 
It is clear that a point $\Aa$ of minimal period $k$ induces a $k$-cycle;
we will sometimes refer to such an $\Aa$ as a {\it primitive} $k$-periodic point.
The zeros of the polynomial
\[
\Phi_k(x):= \prod_{d\mid k} (f^d(x)-x)^{\mu(k/d)},
\] 
where $\mu$ is the Moebius function, are periodic points of $f$ of period dividing
$k$.  A proof that $\Phi_k(x)$ is a polynomial can be found in \cite{pmp}.

To state our main result, let $g(x)\in\Z[x]$ and let $y\in \R$. 
Recall that a {\it derangement} of a set $S$ is a permutation of $S$ with no fixed points.
By the Chebotarev density theorem (\cite{lenstra}, \cite{Neuk}), the density
\begin{equation}\label{density}
\lim_{y\to \infty} \frac{|\{p\leq y : p \text{ is prime and } g(x) \text{ has no zeros mod } p\}|}{|\{ p\leq y : p \text{ is prime }\}|}
\end{equation}
is precisely the proportion of derangements in the Galois group 
$\text{Gal} (g(x)/\Q)$, considered as  a permutation group by its action 
on the set of roots of $g(x)$. When $\Phi_k(x)$ is separable, which holds 
in general, then its roots are periodic of minimal period $k$ \cite{pmp}.  If we set $g(x)=\Phi_k(x)$ 
in \eqref{density} when $\Phi_k(x)$ is separable, then the proportion of derangements 
in Gal$(\Phi_k(x)/\Q)$ is precisely the density of primes $p$ for which $f \pmod p$ 
does not have a $k$-cycle. 

The Galois groups Gal$(\Phi_k(x)/\Q)$ have been studied previously 
(\cite{Vivaldi}, \cite{Bousch}, \cite{pmp}, \cite{Lau}) and in particular
Morton studied these groups in \cite{mortonmain},
proving under strong conditions the isomorphism Gal($\Phi_k(x)/\Q)$ $\cong \Wrr$.
On the other hand, it is not too hard
to show that there exists an  injection of Gal$(\Phi_k(x)/\Q)$ into the wreath product
$\Wrr$ (\cite{Silv} Section 3.9), where $n$ is the number of $k$-cycles of $f$ over $\C$.  
We prove the following results. In particular, Theorem \ref{original}
hold with weaker hypotheses than those in the main theorem of Morton \cite{mortonmain}:

\begin{theorem}\label{newthm}
Let $f(x)$ be a polynomial with integer coefficients.
If $q$ is a power of $p$ and $f(x)\pmod p=x^q+x$,
then $\Phi_{p^n}(x)$ is Eisenstein over $\Q$ for every $n\geq 2$.
\end{theorem}

\begin{theorem}\label{original}
Let $f(x)$ be a polynomial with integer coefficients.
If $\Phi_k(x)$ is separable with $r$ irreducible factors,
then the proportion of primes $p$ such that $f(x) \pmod p$ has
a $k$-cycle is bounded above by $r/k$ for $k>r$.
\end{theorem}

\begin{corollary}\label{maincor}
Let $f(x)$ be as in Theorem \ref{newthm}.  Then
the density of primes $\ell$ such that $f(x) \pmod{\ell}$ has a $p^n$-cycle is
bounded above by $1/p^n$.
\end{corollary}

As an example, the polynomial $f(x)=x^2+3x+2$ has an $8$-cycle modulo $p=389$:
\[
170\longmapsto 237 \longmapsto 88 \longmapsto 230 \longmapsto 299 \longmapsto
52 \longmapsto 139 \longmapsto 290 \longmapsto 170. 
\]
However, there are no primes $p< 389$ such that $f(x)$ has an 8-cycle modulo $p$.
On the other hand of the 3245 primes less than 30,000, $f(x)$
has an 8-cycle modulo 386 of them, or a proportion of about 0.119. For a given prime $p$, it can be hard 
to know whether there is a point of minimal period $k$ for $f(x) \pmod p$.
However, our results show that the density of primes such that $f(x) \pmod p$
has an 8-cycle is at most 1/8.

In Section 2 we provide an overview of the work of Morton {\it et al}.  
We also show the effectiveness (and limitations) of Morton's theorem by example.

In Section 3 we use combinatorics and group theory to study derangements
in subgroups of the generalized symmetric group $\Z/k\Z \text{ wr } S_n$, with
a particular interest in transitive subgroups, or more generally subgroups with 
a fixed number of orbits.  Since $\text{Gal}(\Phi_k(x)/\Q)$
is a subgroup of $\Z/k\Z \text{ wr } S_n$, Theorem \ref{bound} combined
with the Chebotarev density theorem (as we discuss above) yields Theorem \ref{original}.

In Section 4 we show that there exist infinitely many polynomials of given degree with
irreducible $\Phi_k$ for infinitely many $k$ by proving Theorem \ref{newthm}.  In
particular, this combined with our work on derangements in subgroups of 
$\Wrr$ in Section 3 gives us Corollary \ref{maincor}.

\section{Background}

Let $f(x)$ be a polynomial of degree at least 2 with coefficients in a field $K$.
In \cite{pmp} it is shown that $\Phi_k(x) \in K[x]$.  We give some other important
properties of $\Phi_k(x)$ that can be found in \cite{pmp}.

\begin{theorem}
\begin{enumerate}
\item If char $ K \nmid k$, the formula
\begin{equation*}
f^k(x)-x=\prod_{d \mid k} \Phi_d(x)
\end{equation*}
gives a factorization of $f^k(x)-x$ in $K[x]$.  $\deg \Phi_k(x)=\sum_{d\mid k} \mu(k/d)(\deg f)^d$.
\item If char $K \nmid k$, $\Aa$ is a primitive $m$-periodic point of $f$ for $m<k$,
and $\Phi_k(\Aa)=0$, then $(x-\Aa)^2\mid \Phi_k(x)$.
\end{enumerate}
\end{theorem}

Therefore if a non-primitive $k$-periodic point is a root of $\Phi_k(x)$, then it is a multiple root.

Let $\Sigma$ denote
the splitting field of $\Phi_k(x)$ over $K$.  As shown in \cite{Silv} (Section 3.9),
Gal$(\Sigma/K)$ is isomorphic to a subgroup of the wreath product $\Wrr$.

Morton's theorem provides a condition that ensures the Galois group of the 
polynomial $\Phi_k$ attached to certain $f$ is all of $\Wrr$.  For completeness,
we state the theorem over $\Q$, which is enough for the purposes of this paper.
Before stating the theorem, we give an example:

\begin{example}\label{Galois}
Suppose that $f(x)=x^2+1$ and $k=3$.  We calculate the Galois group 
associated to $\Phi_3$.

The primitive 3-periodic points of $f$ satisfy the equation 
\[
\Phi_3(x)=\frac{f^3(x)-x}{f(x)-x}=x^6+x^5+4x^4+3x^3+7x^2+4x+5=0.
\]
Reducing $\Phi_3(x)$ modulo 7 reveals that $\Phi_3(x)$ is irreducible in the polynomial ring
$\Z[x]$.  In fact, the roots of $\Phi_3(x)$ are pairwise distinct and consist of 3 pairs of 
complex conjugates. 
If $K$ denotes the splitting field of $g(x)$ over $\Q$, then $|$Gal$(K/\Q)|\geq 6$,
and since Gal$(K/\Q)\hookrightarrow  \Z/3\Z \text{ wr } S_2$, we have
\begin{equation}\label{stuff}
 |\text{Gal}(K/\Q)|\leq 18.
\end{equation}

To prove that equality holds in \eqref{stuff}, first note that
$g(x) \pmod 5= x(x+4)(x+3)(x^3+4x^2+4x)$, so
Gal$(K/\Q)$ contains a 3-cycle (viewed as a subgroup of $S_6$).  Since
$g(x) \pmod 7$ is irreducible, Gal$(K/\Q)$ contains a 6-cycle, say $\sigma$.
The powers of $\sigma$ never collapse to a 3-cycle so that $|$Gal$(K/\Q)|> 6$.
Hence, $|$Gal$(K/\Q)|=18$.

\end{example}

In order to state Morton's result over $\Q$, we first need some terminology \cite{mortonmain}.  
Let $f(x)$ be a 
monic polynomial in the polynomial ring $\Z[x]$.
Then the discriminant of $\Phi_k(x)$ takes the form
\[
\text{disc } \Phi_k(x)=\pm \Delta^k_{k,k} \cdot \prod_{d \mid k} \Delta^{k-d}_{k,d},
\]
where $\Delta_{k,d}$ satisfy  

\begin{equation}\label{tag'}
\text{Res}(\Phi_k(x),\Phi_d(x))=\pm \Delta_{k,d}^d, \;\; \text{ if } d\mid k, d<k.
\end{equation}
Here, Res$(\cdot,\cdot)$ denotes the resultant of two polynomials.

Morton's theorem over $\Q$ is the following:

\begin{theorem}
Let $f(x) \in \Z[x]$ be monic in $x$ and let $\Sigma$ be the splitting 
field of $f(x)$ over $\Q$.  Assume that
\begin{enumerate}
\item[(i)] $\Phi_k(x)$ is irreducible over $\Q$;
\item[(ii)] Some prime $p$ divides $\Delta_{k,1}$ to the first power but is relatively
prime to $\Delta_{k,d}$ for all $d \neq 1$ and $d \mid k$;
\item[(iii)] $\Delta_{k,k}$ is a square-free integer which is relatively prime
to $\Delta_{k,d}$ for all $d \neq k$ and $d \mid k$, and $\Delta_{k,k}\neq 1$.
\end{enumerate}
Let $\Aa$ be a root of $\Phi_k(x)$ in a splitting field $\Sigma$ of $\Phi_k(x)$ over $\Q$ and
let $kn=$deg $(\Phi_k(x))$.  Then 
the Galois group Gal$(\Phi_k(x)/\Q)=$Gal$(\Sigma/\Q)\cong \Z/k\Z \text{ wr } S_n$.
\end{theorem}

To exemplify just how strong these conditions are, take $f(x)=x^2+3x+2$.
First note, though, that with $k=2$, Morton's theorem is of no real use, 
since one checks that $\Phi_2(x)=
x^2+4x+6$, which is irreducible; hence, Gal($\Phi_2(x)/\Q) \cong \Z/2\Z$.
For periods $k=3,4,5,6,7,8$, for example,
Morton's theorem fails.  Using \eqref{tag'} one easily checks that, in
each case, {\it at least} condition $(iii)$ fails.  For instance,
\begin{eqnarray}\label{tag''}
\text{disc } \Phi_6(x)&=&2^{84}\cdot7^4\cdot 13^{27}\cdot 31^3\cdot41^6\cdot 212901853^6\\
                                    &=&\pm \Delta_{6,6}^6\cdot\Delta_{6,1}^5\cdot\Delta_{6,2}^4\cdot\Delta_{6,3}^3.
\end{eqnarray}
From \eqref{tag'}  one checks that $\Delta_{6,1}=\pm 13$, $\Delta_{6,2}=\pm 7\cdot 13$,
and $\Delta_{6,3}=\pm 2^4\cdot 31$.  In view of \eqref{tag''}, it is now clear that $\Delta_{6,6}$
is not square free (it is divisible by $2^{12}$).

\section{Derangements in Subgroups of Wreath Products}

 In this section our goal is to prove Theorem \ref{original}.

We denote by $d_n$ the number of derangements in the symmetric group on 
$n$ letters, $S_n$.  It is a well-known result that 
\[
\lim_{n\to \infty} \frac{d_n}{n!}=e^{-1}.
\]

We extend this result to the generalized symmetric group $G=\Wrr$.
To begin, write an element in $G$ as $g=((c_1,\dots, c_n), \pi)$, where
$c_i \in \Z/k\Z$, $\pi \in S_n$.  We view $G$ by its action
on $n$ $k$-cycles of some polynomial $f$; suppose
the cycles are $\Oo(\Aa_1),\dots,\Oo(\Aa_n)$, where $\Oo(\Aa_i)=\{\Aa_i,f(\Aa_i), \dots, f^{k-1}(\Aa_i)\}$.
We let $A_i:=\{g\in G : g\Aa_i=\Aa_i\}$ be the pointwise stabalizer of 
the $i$th orbit $\Oo(\Aa_i)$ (if $g$ fixes $\Aa_i$, then $g$ fixes the whole orbit).
If $\A=\A(G)$ denotes the set of derangements in $G$, then 
\[
G\setminus \A=\bigcup_{i=1}^nA_i.
\]
By the principle of inclusion-exclusion,
\[
|\A|=\sum_{I\subseteq\{1,\dots,n\}}(-1)^{|I|}|A_{I}|,
\]
where $A_I:=\bigcap_{i\in I}A_i$.  It is easy to see that for $1\leq \ell \leq n$
\[
|A_{i_1}\cap \cdots \cap A_{i_\ell}|=\binom{n}{\ell}(n-\ell)!k^{n-\ell}= \frac{n!}{\ell!}k^{n-\ell}.
\]
Indeed, there are $\binom{n}{\ell}$ ways to pick an index set of size $\ell$,
$(n-\ell)!$ ways to permute the $n-\ell$ $k$-cycles that are not fixed by $g\in G$
, and finally there are $k$ different orientations for each of the $n-\ell$ $k$-cycles
not fixed by $g\in G$.

We now have
\begin{eqnarray*}
|\A| &=& \sum_{I\subseteq \{1,\dots, n\}} (-1)^{|I|}|A_I|\\
      &=& \sum_{\ell=0}^n (-1)^{\ell}\frac{n!}{\ell!}k^{n-\ell},
      \end{eqnarray*}
from which it follows that  
\begin{equation}\label{tag}
\frac{|\A|}{|G|}=\frac{1}{1!}-\frac{1}{1!k} +\frac{1}{2!k^2}-\cdots+ \frac{(-1)^n}{n!k^n}.
\end{equation}
Taking the limit as $n \to \infty$ gives $e^{-1/k}$.

If $G$ is a permutation group and $\Pp(G):=\frac{|\A(G)|}{|G|}$, then from \eqref{tag} we have
$\Pp(G)=\frac{1}{1!}-\frac{1}{1!k} +\frac{1}{2!k^2}-\cdots+ \frac{(-1)^n}{n!k^n}$ for $G=\Wrr$.
It follows that 
\[
\Pp(G)\geq \frac{k-1}{k} \;\; \text{ for } G=\Wrr.
\]
Since $\Pp(G)=\frac{k-1}{k}+O\left(\frac{1}{k^2}\right)$ as $k\to \infty$,
this inequality is sharp.  In this section, we prove that this
inequality holds for all transitive subgroups of $\Wrr$, and derive an analogous inequality
for subgroups with any number of orbits.

If $G$ is any permutation group, let $m_i(G):=|\{g \in G : ch(g)=i\}|$, where $ch(g)$ is
the permutation character of $g$, defined to be the number of letters fixed by $g$ . 
When $G$ is clear, we will simply write $m_i$.   Define the generating function
 \[
 p_G(t):=\frac{1}{|G|}\sum_{i=0}^{deg(G)}m_it^i.
 \]
\begin{example}\label{dihedral}
Let $G=D_n$, the dihedral group of order $2n$, with its usual action on $n$ letters
visualized as the vertices of a regular $n$-gon.
Then $p_G$ depends on the parity of $n$.  If $n$ is even, then $|\A|= n-1+n/2$, since $n/2$
elements fix two letters, and the identity fixes all letters.  The $n/2$ elements that fix
two letters are the reflections through the $n/2$-lines of symmetry through two vertices
of the regular $n$-gon.  Thus $p_G(t)=\frac{1}{2n}(
\frac{3n-2}{2}+\frac{n}{2}t^2+t^n)$.  If $n$ is odd, then $|\A|=n-1$, $n$ elements fix 1
letter, and the identity fixes all letters.  Thus $p_G(t)=\frac{1}{2n}(n-1+n+t^n)$.
\end{example}

We now give two results from Boston {\it et al}. \cite{Boston},
outlining the basic properties of the polynomial $p_G$:
\begin{theorem} \label{swap}
For any permutation groups $H$ and $K$, we have $\Pp(H \text{ wr } K)=
p_K(\Pp(H))$, where $\Pp(G)=|\A(G)|/|G|$.
\end{theorem}

\begin{theorem} \label{property thm} Suppose $G$ is any permutation group acting on a set $X$.
\begin{enumerate}
\item[(1)]  $p_G(t)$ has degree $n=deg(G)$, and its coefficients are rational, with
denominators dividing $|G|$.
\item[(2)]  $p_G(t)$ and its first $n-1$ derivatives are strictly increasing functions of $t$
for $t\geq 0$.
\item[(3)]  $p_G(0)=\Pp(G)$, $p_G(1)=1$, and for any $x\in X$, $p'_G(1)$ is the number of orbits
of $G$ on $X$.
\item[(4)]  If $G$ is transitive and $H=Stab(x)$ acting on $X-\{x\}$, then $p_H(t)=
p'_G(t)$.
\item[(5)]  The action of $G$ on $X$ is $k$-transitive if and only if $p_G^{(k)}(1)=1$.
\item[(6)]  If $H$ acts on $X$ and $K$ acts on $Y$ then the usual action of $H \times K$ on 
$X\coprod Y$ yields $p_{H\times K}(t)=p_H(t)p_K(t)$.
\item[(7)]  $p_{(H \text{ wr } K)}(t)=p_K(p_H(t))$
\end{enumerate}
\end{theorem}

\begin{example}
If $G=\Z/k\Z \text{ wr }  D_n$ for $n$ even, then by (7) of Theorem \ref{property thm}
and by our remarks in Example \ref{dihedral},  $\Pp(G)=p_{D_n}(\frac{k-1}{k})=
\frac{3n-2}{4n}+\frac{1}{4}\left(\frac{k-1}{k}\right)^2+\frac{1}{2n}\left(\frac{k-1}{k}\right)^n$.
Thus $\lim_{k\to \infty} \lim_{n\to \infty} \Pp(G)=1$.
\end{example}

\begin{example}
Let $q$ be a prime power and let $F_q$ denote the 1-dimensional affine general
linear group over $GF(q)$.  Here $GF(q)$ denotes the finite field with $q$ elements.  
$F_q$ acts on $GF(q)$ as the group of functions $x\mapsto
ax+b$ where $a,b \in GF(q)$, $a\neq 0$.  If $G=\Z/k\Z \text{ wr } F_q$, then 
$\Pp(G)=p_{F_q}(\frac{k-1}{k})$.  To calculate $\Pp(G)$, we note that whenever
$a=1$ and $b\neq 0$, the action $x\mapsto x+b$ is never fixed (i.e. $x+b\neq x$).
We also observe that given $a\neq 0$ or 1, then the action $x\mapsto ax+b$ is 
fixed if and only if $x=\frac{b}{1-a}$.  In the first case, we count $q-1$ such
maps, and in the second case we count $q(q-2)$ of them, leaving
1 element, the identity, which fixes every $x\in GF(q)$.  Thus,
$p_{F_q}(t)=\frac{1}{q(q-1)}(q-1+q(q-2)t+t^q)$, and so
\[
p_{F_q}\left(\frac{k-1}{k}\right)=\frac{1}{q}+\frac{q-2}{q-1}\frac{k-1}{k}+\frac{(\frac{k-1}{k})^q}
{q(q-1)}.
\]
Evidently, $\lim_{q\to \infty}p_{F_q}(\frac{k-1}{k})=\frac{k-1}{k} \to 1$ as $k\to \infty$.
\end{example}

Let $G$ be a permutation group acting on a set $X$.
Recall Burnside's Lemma, which says that the average number of fixed points in
$X$ of an element of $G$ equals the number of orbits of $G$, say $r$:
\[
\frac{1}{|G|}\sum_{g\in G}ch(g)=r.
\]
If $G$ is a subgroup of the generalized symmetric group $\Wrr$, then
elements of $G$ fix $k$ points in $X$ at a time.  If $k$ is large,
this means a large proportion of elements in $G$ cannot fix
any $x$ in $X$.  We make this precise in the following
theorem.

\begin{theorem}\label{bound}
Let $\Gg$ be a subgroup of $\Wrr$, and suppose that $\Gg$ has
$r$ orbits.  Then
there is always an inequality
\begin{equation}\label{tag'''}
\Pp(\Gg)\geq \frac{k-r}{k}
\end{equation}
on the proportion of derangements in $\Gg$.  Moreover, this inequality
is sharp. 

\end{theorem}
 
 \begin{Proof}
 Let $\Gg$ be a subgroup of the wreath product $\Z/k\Z \text{ wr } S_n$.
Denote by $\Qq(\Gg)$ the proportion of non-derangements in
 $\Gg$, and  let $x$ be a randomly chosen
element in $\Gg$.  Let $T_i := \{ x\in \Gg : ch(x)=i \}$, and let us denote by $m_i$ 
the number $|T_i|$. 
Since the $T_i$ partition $\Gg$, we have $\Qq(\Gg)= \frac{1}{|\Gg|} \sum_{i=k,2k,\dots,nk} m_i$  

If $r$ denotes the number of orbits of $\Gg$, then
\begin{eqnarray*}
r=p'_{\Gg}(1) &=& \frac{1}{|\Gg|}\sum_{i=k,2k,\dots,nk} im_i\\
                         &=& \frac{k}{|\Gg|}(m_k+2m_{2k}+\cdots+nm_{nk})\\
                         &>& k\Qq(\Gg),
\end{eqnarray*}
hence $\Qq(\Gg)<\frac{r}{k}$.  Since $\Pp(\Gg)+\Qq(\Gg)=1$, \eqref{tag'''} follows.

It remains to show that \eqref{tag'''} is sharp.  If $A$ is a subgroup of $S_n$ with 
$r$ orbits, then $\Gg':=\Z/k\Z \text{ wr } A$ is a subgroup of $\Wrr$ with $r$
orbits.  As was noted at the beginning of Section 3.2, $\Pp(\Wrr)=\frac{k-1}{k}+
O(\frac{1}{k^2})$. Therefore, $\Pp(\Gg')=\big(\frac{k-1}{k}+O(\frac{1}{k^2})\big)^{r}$.
By the binomial theorem, we have
\[
\left(\frac{k-1}{k}+O\left(\frac{1}{k^2}\right)\right)^{r}= \left(\frac{k-1}{k}\right)^{r}+
O\left(\frac{1}{k^2}\right).
\]
Another application of the binomial theorem gives
\begin{eqnarray*}
\left(\frac{k-1}{k}\right)^{r}&=&\sum_{m=0}^{r}\binom{r}{m}(-1)^mk^{-m}\\
                                      &=& 1-\frac{r}{k}+\frac{\binom{r}{2}}{k^2}-\cdots+\frac{(-1)^{r}}{k^{r}}\\
                                      &=& \frac{k-r}{k} +O\left(\frac{1}{k^2}\right).
                                      \end{eqnarray*}
Thus $\Pp(\Gg')=\frac{k-r}{k} +O(\frac{1}{k^2})$ as $k\to \infty$, and the result follows.
\end{Proof}

Suppose $\Phi_k(x)$ is separable with $r$ irreducible factors.  Then Gal$(\Phi_k(x)/\Q)$
is a subgroup of $\Wrr$ with $r$ orbits, so it satisfies \eqref{tag'''}.   
Then by the Chebotarev density theorem, the proportion of primes $p$ such
that $f(x)\pmod p$ has a $k$-cycle is bounded above by $r/k$.  This proves
Theorem \ref{original}.

\section{Polynomials with Irreducible $\Phi_k$ for
Infinitely Many $k$}

In the previous section we proved that whenever $\Phi_{k}(x)$ is irreducible,
there is a bound on the density of primes $p$ such that $f(x) \pmod p$
has a $k$-cycle.  Namely, we showed that
\[
\Qq(\Gg)\leq  \frac{1}{k},
\]
where $\Gg= \text{Gal}(\Phi_k(x)/\Q)$
and $\Qq(\Gg)$ denotes the proportion of non-derangements in $\Gg$.

In this section we prove Theorem \ref{newthm}.  To do so, we begin with a lemma from
\cite{pmp}, whose proof will be useful for our proof of Theorem \ref{newthm}.

\begin{lemma}\label{fundamental}
If $q=p^s$ for some prime $p$ and $f(x)=x^q+x$, then $\Phi_k(x)$ is 
Eisenstein over $\Q$ whenever $k$ is a power of $p$.
\end{lemma}

\begin{Proof}
Fix the prime $p$ in $f(x)=x^q+x$; we show that $\Phi_k(x)$ is 
Eisenstein with respect to $p$.  The operators $\pi(x):=x^q$ and
$1(x):=x$ are commutative modulo $p$, so by the binomial 
theorem
\begin{eqnarray*}
f^n(x) &\equiv& (\pi+1)^n(x)\\
            &\equiv& \sum_{m=0}^n \binom{n}{m}\pi^m(x)\\
            &\equiv& \sum_{m=0}^n\binom{n}{m}x^{q^m} \pmod p
\end{eqnarray*}
  Thus if $n=p^s$, we have
\begin{equation}
f^n(x) \equiv x^{q^{p(s)}}+x \pmod p
\end{equation}
where $p(s):=p^s$, so 
\begin{equation}
\Phi_n(x)=\frac{f^{p^s}(x)-x}{f^{p^{s-1}}(x)-x}\equiv x^{q^{p(s)}-q^{p(s-1)}} \pmod p.
\end{equation}

It remains to show that $p$ exactly divides $\Phi_n(0)$.  To this end, we prove the congruence
\[
f^n(x) \equiv nx^q +x \pmod{x^{q+1}}.
\]
This is true for $n=1$, so assume it holds for $n=k$.  Write the congruence as
\[
f^k(x)=kx^q+x+x^{q+1}g(x)
\]
for some $g(x)$ in $\Z[x]$.  Then $f^{k+1}(x)=k(x^q+x)^q+x^q+x+f(x)^{q+1}g(f(x))\equiv
(k+1)x^q+x \pmod{x^{q+1}}$, and the result follows by induction.  We now have
\[
\Phi_{p^s}(x)=\frac{p^sx^q+x^{q+1}g_1(x)}{p^{s-1}x^q+x^{q+1}g_2(x)}=\frac{p^s+xg_1(x)}{
p^{s-1}+xg_2(x)}.
\]
Hence $\Phi_{p^s}(0)=p$, and the theorem is proved.
\end{Proof}

We now extend this result to an infinite family of polynomials, which contains
$f(x)=x^q+x$.

\begin{theorem}
Fix a prime $p$ and let $q$ be a power of $p$.
For every polynomial in the infinite family $\F := \{ f(x) \in \Z[x] : f(x) \pmod p=x^{q}+x\}$
the corresponding polynomial $\Phi_k(x)$ is Eisenstein over $\Q$ whenever $k=p^s$, 
$s\geq 2$.  
\end{theorem}

\begin{Proof}
By Lemma \ref{fundamental}, $\Phi_{p^s}(x) \equiv x^{q^{p(s)}-q^{p(s-1)}} \pmod{p}$, so we need
only show that $p$ exactly divides $\Phi_{p^s}(x)$ for every $s\geq 2$. We show that if $p^r$ 
exactly divides $f^{p^s}(0)$, then $p^{r+1}$ exactly divides $f^{p^{s+1}}(0)$, which suffices
since $\Phi_{p^s}(x)\in \Z[x]$,
\[
\Phi_{p^s}(0)=\frac{f^{p^s}(0)}{f^{p^{s-1}}(0)},
\]
and $p\mid f^p(0)$.

Let $p^r$ be the 
largest power of $p$ that divides $f^{p^{s-1}}(0)$, and write $f^{p^{s-1}}(0)=p^rC_1$
for some integer $C_1$ such that $p \nmid C_1$.
In the Taylor expansion
\[
f^{(p-1)\cdot p^{s-1}}(x)= \sum_{\ell \geq 0}\frac{(f^{(p-1)\cdot p^{s-1}})^{(\ell)}(0)}{\ell !}x^{\ell},
\]
we set $x=f^{p^{s-1}}(0)$, giving
\begin{eqnarray*}
f^{p^s}(0) &=& \sum_{\ell \geq 0} \frac{(f^{(p-1)p^{s-1}})^{(\ell)}(0)}{\ell !}(f^{p^{s-1}}(0))^{\ell}\\
                  &=& f^{(p-1)p^{s-1}}(0)+(f^{(p-1)p^{s-1}})'(0)f^{p^{s-1}}(0)+\frac{1}{2}
                          (f^{(p-1)p^{s-1}})''(0)(f^{p^{s-1}}(0))^2+\cdots.
\end{eqnarray*}
Since $f$ is a polynomial, one sees that $\ell ! \mid (f^n)^{(\ell)}(0)$.  By hypothesis,
the polynomial $f(x)$ takes the form
\[
f(x)=x^{q}+a_{q-1}x^{q-1}+\cdots + a_{1}x +a_{0}
\]
where $a_i \equiv 0 \pmod p$ for every $i\neq 1$, while
$a_1 \equiv 1 \pmod p$.  Let $a_1=1+pa$ and $a_0=pb$.
By induction on $i$, one easily obtains the congruence
$f^i(0)\equiv pb(pa+i) \pmod{p^2}$.  Thus with $i=p$ we see that
$f^{p}(0)$ is at least divisible by $p^2$.  It follows that 
all terms following the first two summands are divisible by at least $p^{2r} \geq p^{r+2}$. 
 
We now analyze the first two summands in the Taylor expansion, starting with the
second.  Our first claim is
that the derivative $(f^{jp^{s-1}})'(0) \equiv 1 \pmod{p^2}$ for every $1\leq j \leq p-1$.
We have
\begin{eqnarray*}
(f^{jp^{s-1}})'(0)&=& \prod_{i=0}^{jp^{s-1}-1}f'(f^i(0))\\
                          &=& \prod_{i=0}^{jp^{s-1}-1}(q(f^i(0))^{q-1}+\cdots +2a_2f^i(0)+a_{1})\\
                          &\equiv& (a_1)^{jp^{s-1}} \pmod{p^2},
                       \end{eqnarray*}
where the last equivalence follows from the fact that $p \mid f^i(0)$ for every $i\geq 0$
and $p \mid a_i$ for every $i\neq 1$.
But now note that $(a_{1})^{jp^{s-1}}=(a_{1}^{p^{s-1}})^j \equiv 1 \pmod{p^2}$,
since if $x\equiv 1 \pmod p$ then $x^p \equiv 1 \pmod{p^2}$.  This fact can be seen by
writing $x=1+pn$ and expanding $(1+pn)^p$, noting that $p \mid \binom{p}{s}$ for
$2\leq s \leq p-1$.  

As for the first summand in the Taylor expansion, we prove that
\[
f^{jp^{s-1}}(0) \equiv jp^rC_1 \pmod{p^{r+2}}
\]
for every $1 \leq j \leq p-1$ by induction on $j$.  The result certainly holds
for $j=1$, therefore assume it holds for $j-1 \leq p-1$.  We have
\begin{eqnarray*}
f^{jp^{s-1}}(0) &=& f^{(j-1)p^{s-1}}(f^{p^{s-1}}(0))\\
                          &=& f^{(j-1)p^{s-1}}(0)+(f^{(j-1)p^{s-1}})'(0)f^{p^{s-1}}(0)+\frac{1}{2}(f^{(j-1)p^{s-1}})''(0)
                                   (f^{p^{s-1}}(0))^2+\cdots  \\
                          &=& (j-1)p^rC_1+p^{r+2}C_2 +(p^2n+1)p^rC_1+p^{2r}C_3   \qquad(C_2,C_3 \text{ nonnegative integers})\\
                          &=& jp^rC_1+p^{r+2}C_2+p^{r+2}nC_1+p^{2r}C_3\\
                          &\equiv& jp^rC_1 \pmod{p^{r+2}},
\end{eqnarray*}
as was to be shown.  Now put $j=p$ in
the congruence above: $f^{p^s}(0)=p^{r+1}C_1+p^{r+2}C_2+p^{r+2}nC_1+p^{2r}C_3=
p^{r+1}(C_1+pC_2+pnC_1+p^{r-1}C_3)$, which is exactly divisible by $p^{r+1}$
since $p\nmid C_1$ and $r\geq 2$.
\end{Proof}

It is interesting to note that this theorem fails for $s=1$.  Indeed, take 
$f(x)=x^2+3x+2$ with $p=3$.  Then $\Phi_3(x)=x^6+10x^5
+46x^4+122x^3+199x^2+192x+91$ is actually reducible.


\vspace{12pt}

\noindent{{\bf Acknowledgments.}}
This article is an edited version of the author’s Senior Thesis that was completed when he was a student at the
College of the Holy Cross in 2010. In particular the Reference Section has been updated and expanded. The author would like to thank Professor Rafe Jones for his guidance throughout
this work.  The insights he has shared
and the conversations
we have had have been of great benefit in my refining my own thoughts as presented in this article.

\end{document}